\documentclass[11pt]{article}
\usepackage[cp1251]{inputenc} %if the document contains DOS-codes
\usepackage[dvips]{graphicx}
\usepackage[dvips]{epsfig}
\usepackage{amssymb}
\begin{document}

 \vspace{5mm}
\begin{center}
{\Large{\bf Criteria for analyticity of subordinate semigroups }}

\bigskip
{\bf A. R. Mirotin} \footnote{The author was supported in part by
the State Program of Fundamental Research of Republic of Belarus
under the contract number 20061473.}
\end{center}

\begin{center}
{\small {\bf Abstract}}
\end{center}

 {\small Let $\psi$ be a Bernstein function. A.~Carasso and T.~Kato obtained
necessary and sufficient conditions for $\psi$ to have a property
that $\psi(A)$ generates a quasibounded holomorphic semigroup for
every generator $A$ of a bounded $C_0$-semigroup in a Banach
space, in terms of some convolution semigroup  of measures
associated with $\psi$. We give an alternative to Carasso-Kato's
criterium, and derive several sufficient  conditions for $\psi$ to
have the above-mentioned property.


}

\bigskip
{\bf 1. Introduction}

\medskip
\noindent The well known theorem due to Yosida \cite{Yos60} states
that for every generator $A$ of a bounded $C_0$-semigroup on a
Banach space $X$ its fractional power $-(-A)^\alpha, 0<\alpha<1$
is a generator of a holomorphic semigroup on $X$. The present
paper is devoted to some generalizations and analogs of Yosida's
Theorem in terms of so-called Bochner-Phillips calculus
\cite{Boch, Phil52} (see also \cite[Chap. XIII]{F}; \cite{SMZ98,
Schil,  a&A, BBD}). Though the majority of works on
Bochner-Phillips calculus use the class ${\cal B}$ of (positive)
Bernstein functions, we prefer the class ${\cal T}$ of negative
one. The corresponding reformulation of Bochner-Phillips calculus
is trivial in view of the fact that $\phi(x)\in {\cal B}$ if and
only if $-\phi(-s)\in {\cal T}$.

We say that the function $\psi:(-\infty,0]\to (-\infty,0]$ belongs
to the class ${\cal T}$ of  {\it negative  Bernstein functions} if
$\psi\in C^\infty((-\infty,0))\cap C((-\infty,0])$ and its
derivative is absolutely monotonic, i.e. $\psi^{(n)}\geq 0$ for
all $n\in \mathbb{N}$. It is known that in this case $\psi$
extends analytically to the left half-plane $\Pi_-=\{{\rm
Re}z<0\}$, the extension is continuous on   $\{{\rm Rez}\leq 0\}$,
and has the following integral representation

$$\psi(z)=c_0+\int\limits_{\mathbb{R}_+} (e^{z u}-1)u^{-1}d\rho(u),\quad Rez\leq 0
\eqno(1)
$$
\noindent where $c_0=\psi(0)$, the positive  measure $\rho$ on
$\mathbb{R}_+$ is uniquely determined by $\psi$ and
$\int_{[0,1]}d\rho<\infty,
\int_{[1,\infty)}u^{-1}d\rho(u)<\infty$; the integrand in (1) is
defined for $u=0$  to be equal to $z$ .

Moreover, there is a convolution semigroup $(\nu_t)_{t\geq 0}$  of
sub-probability measures on $\mathbb{R}_+$ with the Laplace transform

$$
g_t(z):=e^{t\psi(z)}=\int\limits_{\mathbb{R}_+} e^{zu}d\nu_t(u),\quad Rez\leq 0
\eqno (2)
$$
\noindent
(see \cite{Sch}, \cite[Chap. XIII]{F}).

The class ${\cal T}$ is a cone which is closed with respect to
compositions and pointwise convergence on $(-\infty, 0]$, and
contains a number of important functions, including (up to affine
changes of variable) fractional powers, the logarithm, the inverse
hyperbolic cosine, and polylogarithms $Li_p$ of all orders $p\in
\mathbb {N}$ \cite{Tula}.

For a negative Bernstein function $\psi$ with integral
representation (1) and a generator $A$ of  a bounded
$C_0$-semigroup $T$ on a complex Banach space $X$ the value of
$\psi$ at $A$ for $x\in D(A)$, the domain of $A$, is defined by
the Bochner integral
$$
\psi(A)x = c_0x + \int\limits_{\mathbb{R}_+} (T(u)-I)xu^{-1}d\rho(u).
$$
\noindent
The closure of this operator, which is also denoted by $\psi(A)$, is a
generator of a bounded  $C_0$-semigroup $g_t(A)$ on $X$ (the "subordinate semigroup"), too.
(For the multidimensional version of this calculus see, e.g., \cite{Mir97}, \cite{Mir98}, \cite{a&A}.)

In the following, without loss of generality we shall assume that
$c_0=0$. The corresponding subclass of  ${\cal T}$ will be denoted
by ${\cal T}_0$. We shall denote also by ${\cal
M}^b(\mathbb{R}_+,\mathbb{C}) ({\cal M}(\mathbb{R}_+,
\mathbb{R}_+))$ the space of all bounded complex valued
(respectively positive) measures on $\mathbb{R}_+$, and by
$C_0(\mathbb{R}_+)$ the space of all continuous complex valued
functions on $\mathbb{R}_+$ which vanish at infinity; $X$ stands
for a complex Banach space.

Another result by Yosida \cite{Yos59} asserts that if the bounded
$C_0$-semigroup $T$ with generator $A$ on  $X$ satisfies
$$
T(t)X\subset D(A), t>0,\quad {\rm and} \quad \lim\sup\limits_{t\downarrow 0} (t\|AT(t)\|)<
\infty, \eqno(Y)
$$
\noindent
then for any $\beta>0$, $e^{-\beta t}T(t)$ can be extended to a bounded  holomorphic
semigroup  on $X$.

We shall denote by  ${\cal T}_Y$   the set of all $\psi\in{\cal T}$  such that
$\psi(A)$ generates  a bounded  $C_0$-semigroup   with property (Y) for every generator
$A$ of     a bounded  $C_0$-semigroup in a Banach space. The class ${\cal T}_Y$ is a cone
\cite[Theorem 6]{CK}. Moreover, it is clear that the composition
$\psi_1\circ\psi_2\in {\cal T}_Y$ if $\psi_1\in {\cal T}_Y$, $\psi_2\in
{\cal T}$. But
 the class ${\cal T}_Y$  is not closed with respect to  pointwise
 convergence.

A. Carasso and T. Kato \cite[ Theorem 4]{CK} obtained necessary
and sufficient conditions for a function $\psi$ to be in ${\cal
T}_Y$ in terms of the semigroup $(\nu_t)_{t\geq 0}$. They also
gave two   necessary  conditions in terms of  $\psi$ itself. Y.
Fujita \cite{Fuj}  obtained sufficient conditions for $\psi$ to be
in ${\cal T}_Y$  in terms of  analytical continuation of $\psi$
and regular variation.

We proceed as follows. First we prove the multiplication rule
which connects the Bochner-Phillips and Hille-Phillips calculi and
then derive the alternative to \cite{CK}  necessary and sufficient
conditions for the inclusion $\psi\in {\cal T}_Y$ (see Theorem 2
below; the variant of this theorem with $C_{00}(\mathbb{R}_+)$
instead of $E(\mathbb{R}_+)$ (for the definition of the last class
see below) first appeared in \cite{SMZ02}). Then we deduce two
theorems from this criterium that give sufficient conditions for
$\psi$ to be in ${\cal T}_Y$ in terms of $\psi$. It should ne
noted that the assumptions of Theorem 4 below contain necessary
conditions, obtained by Carasso and Kato (the idea to employ the
Hausdorff-Young inequality in this context belongs to Carasso and
Kato, too). Finally, we give one more condition, that is
sufficient for the inclusion $\psi\in {\cal T}_Y$. Several
examples have been considered.

\bigskip
{\bf 2. The multiplication rule for  the Bochner-Phillips and
Hille-Phillips calculi, and the criterium for $\psi$ to be in ${\cal T}_Y$}

\medskip
\noindent
In \cite[ Chap.XV]{HPh} the functional calculus (the Hille-Phillips
calculus)  of generators of $C_0$-semigroups have been constructed.  In
particular let   $a\in{\cal M}^b(\mathbb{R}_+,\mathbb{C})$ and
$$
g(s)=La(s):= \int\limits_{\mathbb{R}_+} e^{su}da(u)\quad (s\leq 0)
$$
\noindent
be the {\it Laplace transform of} $a$. Then for a
generator $A$ of  a bounded $C_0$-semigroup $T$ on a complex Banach space $X$ the
value of $g$ at $A$  is the bounded operator on $X$ defined by the Bochner
integral
$$
g(A)x = \int\limits_{\mathbb{R}_+} T(u)xda(u), \quad x\in X.
$$

Our Theorem 1 connects  the Bochner-Phillips and Hille-Phillips
calculi. It is a generalization of Lemma 1 in \cite{SMZ02}. But
first we need the following approximation lemma. We shall denote
by $E(\mathbb{R}_+)$ the complex space of exponential polynomials
of the form
$$
p(t)=\sum\limits_{j=1}^n c_je^{s_jt}, \quad c_j\in \mathbb{C}, s_j<0,
$$
\noindent
endowed with $\sup$-norm on $\mathbb{R}_+$.

\medskip
{\bf Lemma 1}. {\it For every bounded function $q\in C^1(\mathbb{R}_+)$
 with bounded derivative there exists a sequence $q_n\in
E(\mathbb{R}_+)$ such that

{\rm 1)} $q_n\to q$, and $q_n'\to q'$ pointwise on $\mathbb{R}_+$;

{\rm 2)} $(q_n)$ and $(q_n')$ are uniformly bounded on $\mathbb{R}_+$.

\medskip
Proof}. Let us pick a sequence $\tilde{q}_n\in C^1(\mathbb{R}_+)$
such that $\tilde{q}_n(t)=q(t)$ for $t\in[0,n]$,
$\tilde{q}_n(t)=0$ for $t\in[n+1,\infty)$, and   $(\tilde{q}_n)$
and $(\tilde{q}_n')$ are uniformly bounded, $|\tilde{q}_n|<C_1$,
$|\tilde{q}_n'|<C_1$. Define $f_n(x)=\tilde{q}_n(-\log x)$  for
$x\in[0,1]\quad(f_n(0)=0)$. Then $f_n\in C^1([0,1]),\quad
|f_n(x)|<C_1$ for $x\in[0,1]$, and $|f_n'(x)|<C_1 x^{-1}$ for
$x\in(0,1]$. It is well known (see, e. g., \cite[Theorem
8.4.1]{DVL})
 that for every  natural $n$ the algebraic polynomial $p_n$ exists such that
$$
|f_n(x)-p_n(x)|<n^{-1},\quad {\rm and} \quad |f_n'(x)-p_n'(x)|<n^{-1},\quad x\in
[0,1].
$$
\noindent Then $|p_n(0)|<n^{-1}, |p_n(x)|<C_1+1,\quad {\rm and}
\quad  |p_n'(x)|<C_1x^{-1}+1\quad {\rm for} \quad x\in (0,1]$.
Since $f_n(x)=q(-\log x)$  for $x\in(0,1]$, and $n>-\log x$, we
have
$$
|q(-\log x)-p_n(x)|<n^{-1}, \quad  x\in(0,1],\quad n>-\log x.
$$
\noindent
Let $q_n(t):=p_n(e^{-t})-p_n(0)$. Then  $q_n\in
E(\mathbb{R}_+), q_n\to q \quad{\rm on} \quad\mathbb{R}_+$,  and $(q_n)$ and $(q_n')$
are uniformly bounded on $\mathbb{R}_+$.  Finally

$$
|q'(-\log x)(-x^{-1})-p_n'(x)|<n^{-1}, \quad  x\in(0,1], n>-\log x.
$$
\noindent
Putting hear $x=e^{-t}$ we have for all natural $n>t\quad (t\in
\mathbb{R}_+)$  that
$|q'(t)-q_n'(t)|<n^{-1}$.  This completes the proof.

\medskip
For measures $a\in{\cal M}^b(\mathbb{R}_+,\mathbb{C}),\quad {\rm and}\quad
\rho\in{\cal M}(\mathbb{R}_+,\mathbb{R}_+)$   let
$$
K(a,\rho)=\sup\limits_{\phi\in
S}\left|\int\limits_{\mathbb{R}_+}\int\limits_{\mathbb{R}_+}
\phi(r)d_r(a(r-u)-a(r))u^{-1}d\rho(u)\right|
$$
\noindent (if the right hand side  exists), where $S$ is the unit
sphere of the space $E(\mathbb{R}_+)$   with respect to
$\sup$-norm on $\mathbb{R}_+$ . Here we assume that $a=0$ on
$(-\infty;0)$. See the proof of Theorem 5 for an estimate for
$K(a,\rho)$ with bounded positive measure $a$, but
$K(\delta,\delta)=\infty$.

%For example $ K(a,\rho)$ does exist  if the measure
%$a$ is bounded and positive, and the function $u\mapsto
%a([0,u))u^{-1}$ belongs to $L^1(\rho)$; see the proof of Theorem
%5.

 {\bf Theorem 1}. {\it Let $g=La, \quad
a\in{\cal M}^b(\mathbb{R}_+,\mathbb{C}),\quad
 and \quad \psi\in {\cal T}_0$ has integral representation {\rm (1)}. If
$K(a,\rho)<\infty$, then

{\rm 1)} the function $h:=\psi g$ has the form $h=Lb$, where $b
\in{\cal M}^b(\mathbb{R}_+,\mathbb{C})$, $\|b\|= K(a,\rho)$;

{\rm 2)} $g(A)X\subset D(\psi(A))$, $h(A)=\psi(A)g(A)$, and $\|h(A)\|\leq
MK(a,\rho)$ for every operator $A$ in a Banach space $X$, which generates  a
bounded $C_0$-semigroup $T$ with $\|T(t)\|\leq M$.

\medskip
Proof}. Let $a(r)$ denotes the distribution function for $a$,
$a(r)=0$ for $r\in (-\infty,0]$. Then for $s<0$
$$
g(s)= \int\limits_{\mathbb{R}_+}
e^{sr}da(r)=(-s)\int\limits_{\mathbb{R}_+} e^{sr}a(r)dr.
$$
\noindent
Thus for $u\geq 0$ and $s<0$ we have
$$
(e^{su}-1)g(s)=(e^{su}-1)(-s)\int\limits_{\mathbb{R}_+}
e^{sr}a(r)dr\hspace{10cm}
$$
$$
 \hspace{2cm}=(-s)\left(\int\limits_{\mathbb{R}_+}
e^{s(r+u)}a(r)dr-\int\limits_{\mathbb{R}_+} e^{sr}a(r)dr\right)
=(-s)\int\limits_{\mathbb{R}_+} e^{sr}(a(r-u)-a(r))dr=Lb^u(s),
$$
\noindent where $b^u(r)= a(r-u)-a(r)$ has bounded variation and is
concentrated on $\mathbb{R}_+$. Therefore for $\psi\in {\cal T}_0$
with integral representation (1) we get
$$
h(s)=\int\limits_{\mathbb{R}_+}(e^{su}-1)g(s)u^{-1}d\rho(u)=\int\limits_{\mathbb{R}_+}
\int\limits_{\mathbb{R}_+} e^{sr}db^u(r)u^{-1}d\rho(u). \eqno(3)
$$

For $\phi\in E(\mathbb{R}_+)$ let
$$
b(\phi):=\int\limits_{\mathbb{R}_+} b^u(\phi)u^{-1}d\rho(u)=
\int\limits_{\mathbb{R}_+}
\int\limits_{\mathbb{R}_+}\phi(r)d_r(a(r-u)-a(r))u^{-1}d\rho(u)
$$
\noindent be the linear functional on  $E(\mathbb{R}_+)$ (we use
the notation $b^u(\phi)$ for $\int\phi db^u$). By the hypothesis
of the theorem $\|b\|=K(a,\rho)<\infty$, and since
$E(\mathbb{R}_+)$ is dense in $C_0(\mathbb{R}_+)$ by
Stone-Weierstrass Theorem, $b$ extends to a measure $b \in{\cal
M}^b(\mathbb{R}_+,\mathbb{C})$. Furthermore,

$$
b=\int\limits_{\mathbb{R}_+} b^u u^{-1}d\rho(u)
$$

\noindent
(the weak integral; ${\cal M}^b(\mathbb{R}_+,\mathbb{C})$ is endowed with
vague topology).

We claim that for every bounded function $q\in C^1(\mathbb{R}_+)$
with bounded derivative the following equality holds (we write
$b(q)$ instead of $\int_{\mathbb{R}_+}qdb$ in the rest of the
proof)

$$
b(q)=\int\limits_{\mathbb{R}_+} b^u(q)u^{-1}d\rho(u). \eqno(4)
$$

\noindent
In fact, let $(q_n)$ be as in Lemma 1,  and
$|q_n|<C$,  $|q_n'|<C$ for some constant $C>0$.
 Putting $p_n(u):=b^u(q_n)$ we
have
$$
p_n(u)= \int\limits_{\mathbb{R}_+}q_n(r)d_r(a(r-u)-a(r)) =
\int\limits_{\mathbb{R}_+}(q_n(r+u)-q_n(r))da(r). \eqno(5)
$$
\noindent Now let $p(u):=b^u(q)$.  Then $p_n(u)\to p(u)\quad
(n\to\infty)$ pointwise by Lebesgue Theorem. We have
$|q_n(r+u)-q_n(r)|\leq Cu$, and $\leq 2C$. If we take
$w(u)=\min\{u,1\}$, then $w\in L^1(u^{-1}d\rho(u))$ and (5)
implies that $|p_n(u)|\leq 2\|a\|w(u)$. Thus by the  Lebesgue
Theorem
$$
\int\limits_{\mathbb{R}_+} p_n(u)u^{-1}d\rho(u)\to
\int\limits_{\mathbb{R}_+} p(u)u^{-1}d\rho(u) (n\to\infty).
$$
\noindent
On the other hand,
$$
\int\limits_{\mathbb{R}_+} p_n(u)u^{-1}d\rho(u)=
\int\limits_{\mathbb{R}_+} b^u(q_n)u^{-1}d\rho(u)=b(q_n)\to b(q)\quad
(n\to\infty).
$$
\noindent Then $b(q)=\int_{\mathbb{R}_+} p(u)u^{-1}d\rho(u)$, i.
e. (4) holds.  In particular, for $q(r)=e^{sr}\quad (s\leq 0)$ (4)
and (3) imply the equality $h=Lb$ which proves the first statement
of the theorem.

To prove the second one, fix a bounded linear functional $f\in X'$, vector
$x\in D(A)$, and let $q(r)=f(T(r)x)$.  Then $q\in C^1(\mathbb{R}_+)$
and $q$ is bounded together with the derivative $q'(r)=f(T(r)Ax)\quad (r\geq
0)$. For such $q$ equation (4) implies that
$$
f\left(\int\limits_{\mathbb{R}_+}Txdb\right)=\int\limits_{\mathbb{R}_+}f(T(r)x)db(r)=
\int\limits_{\mathbb{R}_+}f\left(\int\limits_{\mathbb{R}_+}T(r)xdb^u(r)\right)u^{-1}d\rho(u).
$$
\noindent So by the definition of the weak integral
$$
\int\limits_{\mathbb{R}_+}\left(\int\limits_{\mathbb{R}_+}T(r)xdb^u(r)\right)u^{-1}d\rho(u)=
\int\limits_{\mathbb{R}_+}Txdb.
$$
In addition, the interior integral in the left hand side here
exists in the sense of Bochner, and
$$
\int\limits_{\mathbb{R}_+}T(r)xdb^u(r)=\int\limits_{[u,\infty)}T(r)xd_ra(r-u)-
\int\limits_{\mathbb{R}_+}T(r)xda(r)\hspace{4cm}
$$
$$
\hspace{2cm}=\int\limits_{\mathbb{R}_+}T(r+u)xda(r)-
\int\limits_{\mathbb{R}_+}T(r)xda(r)= (T(u)-I)g(A)x.
$$
\noindent
Therefore for $x\in D(A)$ we have
$$
h(A)x= \int\limits_{\mathbb{R}_+}T(r)xdb(r)=
\int\limits_{\mathbb{R}_+}(T(u)-I)g(A)xu^{-1}d\rho(u)=\psi(A)g(A)x.
$$
\noindent Since the operator $h(A)$ is bounded, and, on the other
hand, the operator $\psi(A)g(A)$ is closed (as the product of a
closed and a bounded operators), the last equality holds for all
$x\in X$. In particular, $g(A)X\subset D(\psi(A))$. Finally
$$
\|h(A)\|\leq \int\limits_{\mathbb{R}_+}\|T(r)\|d|b|(r)\leq
M\|b\|=MK(a,\rho).
$$
The theorem is proved.

\medskip
{\bf Theorem 2}. {\it Let $\psi\in {\cal T}_0$. Then $\psi\in {\cal T}_Y$
if and only if
$$
K(\nu_t,\rho)=O(t^{-1}),\quad t\downarrow 0 \eqno(6)
$$
\noindent holds (see formulas (1) and (2) for the definitions of
$\rho$ and  $\nu_t$).

\medskip
Proof}. Let (6) holds.  Putting $a=\nu_t$ in Theorem 1 we get that
for sufficiently small $t>0$ the function $h_t=\psi g_t$  has the
form $h_t=Lb_t$, where $b_t$ is a bounded  measure on
$\mathbb{R}_+$, $\|b_t\|=K(\nu_t,\rho)$. In addition,
$g_t(A)X\subset D(\psi(A))$ for all $t>0$ ($\psi(A)=$ generator of
the semigroup $g_t(A)$) and

$$
\|h_t(A)\|=\|\psi(A) g_t(A)\| \leq MK(\nu_t,\rho).
$$
\noindent
Now (6) implies (Y) with $g_t(A)$ instead of $T(t)$.

To prove the converse, consider $X=C_0(\mathbb{R}_+)$ with
$\sup$-norm, let $\psi\in {\cal T}_Y$, and let $T$ be the
$C_0$-semigroup of shifts on $X$, $(T(r)x)(v)=x(v+r)$ (in this
concrete situation $A$ is a derivation with appropriate domain).
Then, for each $x\in C^1(\mathbb{R}_+)\cap C_0(\mathbb{R}_+),\quad
t>0$ integration by parts gives
$$
y(v):=g_t(A)x(v)=\int\limits_{\mathbb{R}_+}x(v+r)d\nu_t(r)=
-\int\limits_{\mathbb{R}_+}x'(v+r)\nu_t(r)dr.
$$
\noindent
Therefore
$$
\psi(A)g_t(A)x(v)=\int\limits_{\mathbb{R}_+}(y(v+u)-y(v))u^{-1}d\rho(u)\hspace{6cm}
$$
$$
\hspace{3cm}=\int\limits_{\mathbb{R}_+}\left(-\int\limits_{\mathbb{R}_+}x'(v+u+r)\nu_t(r)dr
+\int\limits_{\mathbb{R}_+}x'(v+r)\nu_t(r)dr\right)u^{-1}d\rho(u).
$$
\noindent
Since $\nu_t$ is concentrated on $\mathbb{R}_+$, we get
$$
\int\limits_{\mathbb{R}_+}x'(v+u+r)\nu_t(r)dr=\int\limits_{\mathbb{R}_+}x'(v+r)\nu_t(r-u)dr,
$$
\noindent
and  thus
$$
\psi(A)g_t(A)x(v)= \int\limits_{\mathbb{R}_+}\left(\int\limits_{\mathbb{R}_+}
(\nu_t(r)-\nu_t(r-u))x'(v+r)dr\right)u^{-1}d\rho(u).
$$
\noindent But integration by parts  gives since
$\nu_t(0)=\nu_t(-u)=0$,
$$
\int\limits_{\mathbb{R}_+}(\nu_t(r)-\nu_t(r-u))x'(v+r)dr=
\int\limits_{\mathbb{R}_+}x(v+r)d_r(\nu_t(r-u)-\nu_t(r)).
$$
\noindent
Finally,  for each $x\in
C^1(\mathbb{R}_+)\cap C_0(\mathbb{R}_+),\quad v\geq 0$
$$
\psi(A)g_t(A)x(v)= \int\limits_{\mathbb{R}_+}
\left(\int\limits_{\mathbb{R}_+}x(v+r)d_r(\nu_t(r-u)-\nu_t(r))\right)u^{-1}d\rho(u).
$$
\noindent
Taking into account that $t\|\psi(A)g_t(A)\|\leq C$ for some $C>0$
and all $t\in (0,1]$  we have for our $x$ with $\|x\|=1$ that
$|\psi(A)g_t(A)x(v)|\leq Ct^{-1}$. So for each $v\geq 0$,  $t\in (0,1]$
$$
 \left| \int\limits_{\mathbb{R}_+}
\left(\int\limits_{\mathbb{R}_+}x(v+r)d_r(\nu_t(r-u)-\nu_t(r))\right)u^{-1}d\rho(u)\right|
\leq  Ct^{-1}.
$$
\noindent

Since $C^1(\mathbb{R}_+)\cap C_0(\mathbb{R}_+)$ is dense in
$C_0(\mathbb{R}_+)$, it follows for $v=0$ that  $K(\nu_t,\rho)=O(t^{-1}),\quad
t\downarrow 0$,
as desired.

\bigskip
{\bf 3. Sufficient conditions for $\psi$ to be in ${\cal T}_Y$ in terms of $\psi$}

\medskip
In the following we shall denote  by ${\cal F}$ the Fourier transform on
$\mathbb{R}$,
$$
{\cal F}f(\lambda)=\frac{1}{\sqrt{
2\pi}}\int\limits_{\mathbb{R}}e^{-i\lambda t}f(t)dt,
$$
\noindent and by  ${\cal F}^{-1}$ the inverse of ${\cal F}$. Let
$$
F_t(\lambda)=e^{t\psi(i\lambda)}\psi(i\lambda)\quad ({\rm Im}\lambda\geq 0, t>0).
$$
\noindent
The restriction $F_t|\mathbb{R}$  will be also denoted by $F_t$.

\medskip
{\bf Theorem 3}. {\it Let $\psi\in {\cal T}_0$. Assume that

{\rm (i)} the derivative $\partial/\partial yF_t(y)$ exists for a.e. $y\in
\mathbb{R}$ and each sufficiently small $t>0$;

{\rm (ii)} for some $p\in (1,2]$ functions $F_t$ and  $\partial/\partial yF_t$
both belong to
$L^p(\mathbb{R})$ for  each sufficiently small  $t>0$;

{\rm (iii)}  ${\cal F}F_t$ is concentrated on  $\mathbb{R}_+$ for
each sufficiently small $t>0$;

{\rm (iv)} $\|F_t\|_p^{1/q}\|\partial/\partial
yF_t\|_p^{1/p}=O(t^{-1}) \quad as \quad t\downarrow 0\quad (p^{-1}+q^{-1}=1)$.

Then   $\psi\in {\cal T}_Y$.

 \medskip
 Proof.}  First we prove that $f_t:= {\cal F}F_t\in L^1(\mathbb{R}_+)$, and
$F_t={\cal F}^{-1}f_t$. Indeed, $f_t\in L^q(\mathbb{R})$, and
${\cal F}(\partial/\partial yF_t)(y)=iyf_t(y)\in L^q(\mathbb{R})$.
By H\"older's inequality $f_t(y)=(iyf_t(y))(iy)^{-1}\in
L^1(\{|y|>1\})$, and so   $f_t\in L^1(\mathbb{R})$. Now by the
Inverse Theorem for the Fourier transform, $F_t(y)={\cal
F}^{-1}f_t(y)$ a.e. $y\in \mathbb{R}$, and by the continuity the
last equality holds for all $y\in \mathbb{R}$. Therefore we have
for the Laplace transform
$$
Lf_t(z)=\int\limits_{\mathbb{R}_+}e^{zr}f_t(r)dr=\sqrt{2\pi}e^{t\psi(z)}\psi(z),\quad
{\rm Re}z\leq 0,
$$
\noindent because  both sides here are analytic on the left
half-plane $\Pi_-$, continuous on its closure, and coincide on its
boundary $i\mathbb{R}$. In particular, $Lf_t(s)=
\sqrt{2\pi}e^{t\psi(s)}\psi(s)$ for all $s\leq 0$. It follows that
for an arbitrary exponential polynomial $\phi\in E(\mathbb{R}_+),
\quad \phi(r)=\sum_j c_je^{s_jr}\quad (c_j\in \mathbb{C}, s_j<0)$
we have
$$
\int\limits_{\mathbb{R}_+}\phi(r)f_t(r)dr =\sqrt{2\pi}\sum\limits_j
c_je^{t\psi(s_j)}\psi(s_j).
$$

On the other hand,
$$
\int\limits_{\mathbb{R}_+}\phi(r)d_r(\nu_t(r-u)-\nu_t(r))=
\int\limits_{[-u,\infty)}\phi(r+u)d\nu_t(r)-\int\limits_{\mathbb{R}_+}\phi(r)d\nu_t(r)
$$
$$
\hspace{7cm}=\int\limits_{\mathbb{R}_+}(\phi(r+u)-\phi(r))d\nu_t(r)=
\sum\limits_j c_j(e^{s_ju}-1)e^{t\psi(s_j)},
$$
\noindent
and thus
$$
\int\limits_{\mathbb{R}_+}\int\limits_{\mathbb{R}_+}\phi(r)d_r(\nu_t(r-u)-\nu_t(r))
u^{-1}d\rho(u)=\sum\limits_j
c_je^{t\psi(s_j)}\psi(s_j).
$$
\noindent
Now we conclude that ($E(\mathbb{R}_+)$ is dense in $C_0(\mathbb{R}_+)$)
$$
K(\nu_t,\rho)=\frac{1}{\sqrt{2\pi}}\sup\limits_{\phi\in
S}\left|\int\limits_{\mathbb{R}_+}\phi(r)f_t(r)dr\right|=\frac{1}{\sqrt{2\pi}}\|f_t\|_1.
$$

Let $k_t(u):=iuf_t(u)$. Then $k_t= {\cal F}(\partial/\partial yF_t)$, and
using the Hausdorff-Young inequality we obtain
$$
 \|f_t\|_q\leq \|F_t\|_p,\quad   \|k_t\|_q\leq \|\partial/\partial
yF_t\|_p.
$$
\noindent Next, for any $v>0$ H\"older's inequality gives
$$
\int\limits_{[0,v]}|f_t(u)|du\leq  \|f_t\|_q v^{1/p},
$$
$$
\int\limits_{[v,\infty)}(u|f_t(u)|)u^{-1}du\leq  \|k_t\|_q
(p-1)^{-1/p}v^{-1/q}.
$$
\noindent
Then for any $v>0$
$$
\|f_t\|_1\leq  \|f_t\|_q v^{1/p}+\|k_t\|_q (p-1)^{-1/p}v^{-1/q} \leq
\|F_t\|_p v^{1/p}+ \|\partial/\partial yF_t\|_p (p-1)^{-1/p}v^{-1/q}.
$$
Therefore, on choosing $v=(p-1)^{1/q} \|\partial/\partial yF_t\|_p /
\|F_t\|_p$, it follows that
$$
K(\nu_t,\rho)=\frac{1}{\sqrt{2\pi}}\|f_t\|_1\leq {\rm const} \|F_t\|_p^{1/q}
\|\partial/\partial  yF_t\|_p^{1/p}=O(t^{-1})\quad {\rm as} \quad t\downarrow 0.
$$
\noindent Application of Theorem 2 completes the proof.

\medskip
Before  formulating the next theorem we note that by
 \cite[Theorem 4]{CK} every $\psi\in {\cal T}_Y\cap {\cal T}_0$ maps $\Pi_-$ into
a  truncated sector
$$
S(\theta,\beta):= (\beta+\{|\arg(-z)|<\theta\})\cap\Pi_-
$$
\noindent for some $\beta\geq 0, \theta\in(0,\pi/2)$, and there
exist constants $k, k>0$, and $\gamma, \gamma\in (0,1)$, such that
$|\psi(z)|\leq k|z|^\gamma,\quad |z|\geq 1,\quad {\rm Re}z\leq 0$.
The problem is what one can add to this conditions to obtain
(necessary and) sufficient conditions for   $\psi$ to be in ${\cal
T}_Y$. Now we shall deduce the partial answer to this question
from Theorem 3.

\medskip
{\bf Theorem 4}. {\it Let $\psi\in {\cal T}_0$, and assume that
the following conditions hold:

{\rm (i)}  $\psi:\Pi_- \to  S(\theta,\beta)$ for some  $\beta\geq 0,
\theta\in(0,\pi/2)$;

there exist such positive constants $k, b, \alpha, \gamma$ and
$R$, that $\alpha\leq \gamma <1$, $R\geq 1$, and

{\rm (ii)}  $b|z|^\alpha\leq |\psi(z)|\leq k|z|^\gamma$ for
$z\in\Pi_-,\quad |z|\geq R$;

{\rm (iii)} the function $y\mapsto \psi(iy)$ is differentiable for a. e.
$y\in \mathbb{R}$ and
$$
|\psi'(iy)|\leq k|y|^\delta,\quad a. e. \quad y\in \mathbb{R}, \quad |y|\geq R,
$$
\noindent
for some $\delta\in(\alpha-\gamma-1, 2\alpha-\gamma-1)$ if
$\alpha<\gamma$, and $\delta=\gamma-1$ if  $\alpha=\gamma$;

{\rm (iv)}  $\psi'(iy)\in L^p([0,R])$  for some $p\in(1,2]$ such that
$p=\min\{2,(\alpha-\gamma-\delta)^{-1},(\alpha-\delta-1)/(\gamma-\alpha)\}$
if $\alpha<\gamma$, and $p<\min\{2, (1-\gamma)^{-1}\}$ if  $\alpha=\gamma$.

Then $\psi\in {\cal T}_Y$.

\medskip
Proof}. We shall verify  all the conditions of Theorem 3 for
$\psi$.  Let $a_1=\max\{|\psi(z)|| z\in\Pi_-, |z|\leq R\}$,
$m_1=\min\{|\psi(z)|-b|z|^\alpha| z\in\Pi_-, |z|\leq R\}$. Then
$b|z|^\alpha +a_2\leq |\psi(z)|\leq k|z|^\gamma +a_1$ for
$z\in\Pi_-$, where $a_2=\min\{0,m_1\}$. Since $\psi(iy)-\beta\in
S(\theta,0)$, we have $-{\rm Re}\psi(iy)+\beta\geq
\cos\theta(|\psi(iy)|-\beta)$, and ${\rm Re}\psi(iy)\leq
-c_1|y|^\alpha +c_2$, where $c_1=b\cos\theta>0, c_2\in
\mathbb{R}$. It follows that
$$
|F_t(y)|\leq e^{c_2t}e^{-c_1t|y|^{\alpha}}(k|y|^\gamma+a_1),
$$
\noindent
and ($p\geq 1$)
$$
\|F_t\|_p\leq  e^{c_2t}2^{1/p}\left(\int\limits_{\mathbb{R}_+}
e^{-c_1pt|y|^{\alpha}}(k|y|^\gamma+a_1)^pdy\right)^{1/p}.
$$
\noindent
Putting $x=ty^\alpha$ we get for some constant $c_3>0$

$$
\|F_t\|_p\leq  c_3e^{c_2t}t^{-\gamma /\alpha-1/{\alpha p}}\left(\int\limits_{\mathbb{R}_+}
e^{-c_1px}(kx^{\gamma /\alpha}+a_1t^{\gamma /\alpha})^px^{1/\alpha -1}dx\right)^{1/p}.
$$
\noindent The integral converges for all $t\geq 0, p\geq 1$, and
by B. Levi's Theorem
$$
\|F_t\|_p=O(1)t^{-\gamma /\alpha - 1/\alpha p} \quad{\rm as}\quad
t\downarrow 0. \eqno(7)
$$
Let    $\alpha<\gamma$,
$p=\min\{2,(\alpha-\gamma-\delta)^{-1},(\alpha-\delta-1)/(\gamma-\alpha)\}$,
$\delta\in(\alpha-\gamma-1, 2\alpha-\gamma-1)$. Since
$$
|\partial/\partial y F_t(y)|\leq  e^{c_2t}  e^{-c_1t|y|^{\alpha}}
(tk|y|^\gamma+ta_1+1)|\psi'(iy)|,
$$
\noindent
we have
$$
\|\partial/\partial y F_t\|_p\leq
e^{c_2t}2^{1/p}\Bigl(\int\limits_{[0,R]}e^{-c_1pt|y|^{\alpha}}
(tk|y|^\gamma+ta_1+1)^p|\psi'(iy)|^pdy
$$
$$
+k^p\int\limits_{[R,\infty)}e^{-c_1pt|y|^{\alpha}}
(tk|y|^\gamma+ta_1+1)^py^{\delta p}dy\Bigr)^{1/p}.
$$
\noindent
Putting $x=ty^\alpha$ in the second integral we get
$$
\|\partial/\partial y F_t\|_p\leq
e^{c_2t}2^{1/p}t^{-\frac{\gamma+\delta}{\alpha} - \frac{1}{\alpha
p}+1}\Bigl(t^{(\frac{\gamma+\delta}{\alpha}+ \frac{1}{\alpha
p}-1)p} \int\limits_{[0,R]}e^{-c_1pt|y|^{\alpha}}
(tk|y|^\gamma+ta_1+1)^p|\psi'(iy)|^pdy
$$
$$
+ k^p\alpha^{-1}\int\limits_{[tR^\alpha,\infty)}e^{-c_1px}
(kx^{\frac{\gamma}{\alpha}}+
t^{\frac{\gamma}{\alpha}-1}(ta_1+1))^p x^ {\frac{\delta
p+1}{\alpha} -1}dx\Bigr)^{1/p} .\eqno(8)
$$
\noindent The second integral in (8) converges for all $t\geq 0$
because $(\gamma+\delta)p /\alpha+1/\alpha-1>-1$  for our $p$ and
$\delta$. Note that $(\gamma+\delta) /\alpha + 1/\alpha p-1\geq
0$. Therefore (8) implies
$$
\|\partial/\partial y F_t\|_p= O(1)t^{-(\gamma+\delta) /\alpha -
1/\alpha p+1} \quad{\rm as}\quad t\downarrow 0. \eqno(9)
$$
\noindent
It follows from (7) and (9)  that   for our $\delta$  we have
$$
\|F_t\|_p^{1/q}\|\partial/\partial
yF_t\|_p^{1/p}= O(1)t^{-\gamma /\alpha -
1/\alpha p-(\delta/\alpha-1)/p} =O(t^{-1}) \quad {\rm as} \quad t\downarrow 0,
$$
\noindent
because  $\gamma /\alpha +1/\alpha p+(\delta/\alpha-1)/p\leq 1$.

The case $\gamma=\alpha, \delta=\gamma-1, 1<p<\min\{2,(1-\gamma)^{-1}\}$
can be examined in the same manner.

Finally since  $\psi(i\lambda)-\beta\in S(\theta,0)$ for
$\lambda\in \mathbb{C}$ with ${\rm Im}\lambda\geq 0$, we have  for
such $\lambda$ (as above)

$$
|F_t(\lambda)|\leq
e^{c_2t}e^{-c_1t|\lambda|^{\alpha}}(k|\lambda|^\gamma+a_1).
$$
\noindent

Then for $t>0$ ($\lambda =s+iy, y>0$)
$$
\int\limits_{\mathbb{R}}|F_t(s+iy)|ds\leq 2e^{c_2t}
\int\limits_{\mathbb{R}_+}
e^{-c_1t(s^2+y^2)^{\alpha/2}}(k(s^2+y^2)^{\gamma/2} +a_1)ds
$$
$$
\hspace{3cm}\stackrel{[s^2+y^2=v]}{=}e^{c_2t}\int\limits_{[y^2,\infty)}
e^{-c_1tv^{\alpha/2}}(kv^{\gamma/2}+a_1)(v-y^2)^{-1/2}dv.
$$
\noindent
But
$$
\int\limits_{[y^2,y^2+1]}
e^{-c_1tv^{\alpha/2}}(kv^{\gamma/2}+a_1)(v-y^2)^{-1/2}dv\leq\max
\limits_{v\geq 0}e^{-c_1tv^{\alpha/2}}(kv^{\gamma/2}+a_1)
\int\limits_{[0,1]}u^{-1/2}du.
$$
\noindent
Furthermore
$$
\int\limits_{[y^2+1,\infty)}
e^{-c_1tv^{\alpha/2}}(kv^{\gamma/2}+a_1)(v-y^2)^{-1/2}dv\leq
 \int\limits_{[1,\infty)}
e^{-c_1tv^{\alpha/2}}(kv^{\gamma/2}+a_1)dv.
$$
\noindent Thus $F_t$ belongs to the Hardy class $H^1(\{{\rm
Im}\lambda
>0\})$ for all $t>0$ and therefore ${\cal F}F_t$ is concentrated
on $\mathbb{R}_+$. This completes the proof.

\medskip
{\bf Corollary 1}. {\it Let $\psi\in {\cal T}_0$, and assume that
the following conditions hold:

{\rm (i)}  $\psi:\Pi_- \to  S(\theta,\beta)$ for some  $\beta\geq 0,
\theta\in(0,\pi/2)$;

{\rm (ii)}  $\psi(z)\asymp z^\gamma$ for some $\gamma\in(0,1)$ ($z\to\infty,
\quad z\in\Pi_-$);

{\rm (iii)} the function $y\mapsto \psi(iy)$ is differentiable for a. e.
$y\in \mathbb{R}$ and
$$
|\psi'(iy)|\leq k|y|^{\gamma -1},\quad a. e. \quad y\in \mathbb{R}.
$$

Then $\psi\in {\cal T}_Y$. }

\medskip
{\it Example 1} \cite{Yos60}. Let
$\psi(z)=c^\alpha-(c-z)^\alpha,\quad \alpha\in(0,1),\quad c\geq
0$. In this case, all the conditions of  Corollary 1 (and hence of
Theorems 3 and 4) are clear.

\medskip
Now we shall give an example of a function  $\psi\in{\cal T}_0$
that satisfies all the conditions of Theorem 4, but conditions of
the Theorem in \cite{Fuj} do not hold for $-\psi(-x)$.

{\it Example 2}.  Let $0<\alpha<\beta<1$, and
$$
\psi(z)=-(-z)^\alpha+(e^{-(-z)^\beta} -1).
$$
\noindent Since the summands map $\Pi_-$ into a sector and into a
truncated sector respectively, the condition (i) of Theorem 4
holds. It is easy to verify that $\psi(z)\sim z^\alpha\quad{\rm
as}\quad z\to\infty, z\in\Pi_-$, $\psi'(iy)\sim\alpha
|y|^{\alpha-1}$ as $y\to\infty$. Finally (iv) holds for $p\in
(1,\min\{2,(1-\alpha)^{-1}\})$.  At the same time, $-\psi(-x)$ is
not regularly varying.

\bigskip
{\bf 4. Further sufficient conditions for $\psi$ to be in ${\cal
T}_Y$}

\medskip
 In this section, we shall deduce  further conditions from
Theorem 2, that are sufficient for  $\psi\in{\cal T}_Y$.

{\bf Theorem 5}. {\it Let $\psi\in{\cal T}_0$ and the function
$r\mapsto \nu_t([r-u,r))$ is monotone decreasing on $[u,+\infty)$
$(u\geq 0)$ for each sufficiently small $t>0$. If
$$
\int\limits_{\mathbb{R}_+}\nu_t([0,u))u^{-1}d\rho(u)=O(t^{-1})
\quad {\rm as}\quad t\downarrow 0,
$$
\noindent
then $\psi\in{\cal T}_Y$.

Proof}. Let $a\in{\cal M}^b(\mathbb{R}_+,\mathbb{R}_+)$, and the
function  $r\mapsto a([r-u,r))$ is monotone decreasing on
$[u,+\infty)$   $(u\geq 0)$. Since $\lim_{r\to
+\infty}a([r-u,r))=0$ for every $u>0$, for all and $\phi\in
E(\mathbb{R}_+)$ with $\sup|\phi|\leq 1$ we find ($a(r)=a([0,r))$
for $r>0$, and $a(r)=0$ for $r\in (-\infty,0]$)
$$
|\int\limits_{\mathbb{R}_+}\phi(r)d_r(a(r-u)-a(r))|\leq
Var_{r\in\mathbb{R}_+}(a(r-u)-a(r))=
$$
$$
Var_{r\in [0,u)}a(r)+Var_{r\in
[u,+\infty)}(a(r-u)-a(r))=2a([0,u)).
$$

Thus
$$
K(a,\rho)\leq 2\int\limits_{\mathbb{R}_+}a([0,u))u^{-1}d\rho(u).
$$
\noindent  It remains to put here $a=\nu_t$ and to apply Theorem
2.

\medskip
{\it Example 3} (cf. \cite[Example 1]{CK}). Let for $b>0$
$$
\psi(z)=\log b-\log(b-z).
$$
\noindent
It is well known that  $d\rho(u)=e^{-bu}du$ and
$$
e^{t\psi(s)}=b^t(b-s)^{-t}=b^t\Gamma(t)^{-1}
\int\limits_{\mathbb{R}_+} e^{sr}r^{t-1}e^{-br}dr.
$$
\noindent So $d\nu_t(r)=b^t\Gamma(t)^{-1}r^{t-1}e^{-br}dr$, and
$\nu_t$ has monotone decreasing density for $t\in (0,1)$.
Therefore
$$
\int\limits_{\mathbb{R}_+}\nu_t([0,u))u^{-1}d\rho(u)=
\int\limits_{\mathbb{R}_+}\left(\int\limits_{[0,u)}b^t\Gamma(t)^{-1}r^{t-1}
e^{-br}dr\right)u^{-1}e^{-bu}du\leq
$$
$$
\hspace{4cm}b^t\Gamma(t)^{-1}
 \int\limits_{\mathbb{R}_+}\left(\int\limits_{[0,u)} r^{t-1}dr\right)u^{-1}e^{-bu}du=
\frac{1}{t}.
$$
Thus $\psi\in{\cal T}_Y$ by Theorem 5.

\medskip
{\it Example 4} (cf. \cite{SMZ02}). Let
$$
\psi(s)={\rm acosh} b-{\rm acosh}(b-s) \quad (b\geq 1, s\leq 0).
$$
\noindent Since $\psi\in{\cal T}_Y$ implies
$-\psi(-c)+\psi(s-c)\in{\cal T}_Y$  for all $c\geq 0$, one can to
restrict ourselves to the case $b=1$.  In this case,
$\psi\in{\cal T}_0$   with
 $d\rho(u)=e^{-u}I_0(u)du$ (the corresponding integral representation (1) can be verified
by differentiation under the integral sign), and $e^{t\psi(s)}=
Lf_t(s)$ with  $f_t(r) = tr^{-1}e^{-r}I_t(r),\quad r>0$ ($I_t$
denotes the Bessel function of the first kind). Hence, $d\nu_t(r)=
f_t(r)dr$, and $\nu_t$ has monotone decreasing density for $t\in
(0,1)$ (see \cite{SMZ02}). The calculations from Example 3 in
\cite{SMZ02} show, that the conditions of Theorem 5 hold.  So,
$\psi\in{\cal T}_Y$.

\bigskip
{\bf Acknowledgement}. The author is sincerely grateful to the
referee for helpful comments and suggestions.

\bigskip
%\label{liter}\def\refname{}

\begin{flushleft}
Department of Mathematics

Skoryna Gomel State University

Sovietskaya, 104

246119  Gomel'

 Belarus'

E-mail: amirotin@yandex.ru
\end{flushleft}

\end{document}